\documentclass[12pt]{article}
\renewcommand{\baselinestretch}{1.5}

\reversemarginpar

\newcommand{\D}{\displaystyle}

\newfam\scfam

\font\sevenrm=cmr7
\font\seveni=cmmi7
\font\sevensy=cmsy7
\font\sevenbf=cmbx7

\font\ninerm=cmr9
\font\ninei=cmmi9
\font\ninesy=cmsy9
\font\ninebf=cmbx9

\font\twelverm=cmr12
\font\twelvei=cmmi12
\font\twelvesy=cmsy10 scaled \magstep1
\font\twelveex=cmex10 scaled \magstep1
\font\twelvebf=cmbx12
\font\twelvesl=cmsl12
\font\twelvett=cmtt12
\font\twelveit=cmti12
\font\twelvecsc=cmcsc10 scaled \magstep1
\ifx\twelvepoint\undefined
   \def\twelvepoint{\def\rm{\fam0\twelverm}% switch to 12-point type
       \textfont0=\twelverm \scriptfont0=\ninerm
\scriptscriptfont0=\sevenrm
       \textfont1=\twelvei  \scriptfont1=\ninei
\scriptscriptfont1=\seveni
       \textfont2=\twelvesy \scriptfont2=\ninesy
\scriptscriptfont2=\sevensy
       \textfont3=\twelveex
\scriptfont3=\twelveex\scriptscriptfont3=\twelveex
       \textfont\itfam=\twelveit  \def\it{\fam\itfam\twelveit}%
       \textfont\slfam=\twelvesl  \def\sl{\fam\slfam\twelvesl}%
       \textfont\ttfam=\twelvett  \def\tt{\fam\ttfam\twelvett}%
       \textfont\bffam=\twelvebf  \scriptfont\bffam=\ninebf
       \scriptscriptfont\bffam=\sevenbf  \def\bf{\fam\bffam\twelvebf}%
       \textfont\scfam=\twelvecsc \def\sc{\fam\scfam\twelvecsc}%
       \normalbaselineskip=14.4pt
       \parindent=24pt
       \setbox\strutbox=\hbox{\vrule height9.5pt depth4.5pt width0pt}%
       \normalbaselines\rm}
   \fi
\catcode`\@=11
\def\undefine#1{\let#1\undefined}
\def\newsymbol#1#2#3#4#5{\let\next@\relax
 \ifnum#2=\@ne\let\next@\msafam@\else
 \ifnum#2=\tw@\let\next@\msbfam@\fi\fi
 \mathchardef#1="#3\next@#4#5}
\def\mathhexbox@#1#2#3{\relax
 \ifmmode\mathpalette{}{\m@th\mathchar"#1#2#3}%
 \else\leavevmode\hbox{$\m@th\mathchar"#1#2#3$}\fi}
\def\hexnumber@#1{\ifcase#1 0\or 1\or 2\or 3\or 4\or 5\or 6\or 7\or 8\or
 9\or A\or B\or C\or D\or E\or F\fi}
    \font\elevenBb=msbm10 scaled 1200
    \font\eightBb=msbm7 scaled 1200 
    \font\sixBb=msbm5 scaled 1200
    \newfam\Bbfam
    \textfont\Bbfam\elevenBb
    \scriptfont\Bbfam\eightBb
    \scriptscriptfont\Bbfam\sixBb
\edef\Bbfam@{\hexnumber@\Bbfam}
    \def\Bb{\fam\Bbfam}
    \def\Q{{\Bb Q}}
    
    \def\vn{\mathchar"0\Bbfam@3F}

\begin{document}
\newtheorem{theorem}{Theorem}
\newtheorem{lemma}{Lemma}
\newtheorem{prop}{Proposition}
\newtheorem{cor}{Corollary}

 \begin{center}
 {\LARGE On the Middle Coefficient of a Cyclotomic Polynomial
 } 
 \vskip0.3in
 
 {\large Gregory P. Dresden}
 \end{center}
 \thispagestyle{empty}
 
 \vskip0.3in

\noindent The cyclotomic polynomials $\Phi_n$ for $n=1,2,3,\dots$
 (familiar to every student of algebra) are
the minimal polynomials for the 
primitive $n${\rm{th}} roots of unity:
\[
\Phi_n(x) = \prod_{(k,n) = 1} \left(x- e^{2\pi ik/n}\right).\nonumber 
\]
Clearly $\Phi_n$ has degree $\phi(n)$, where $\phi$
signifies Euler\rq s totient function. These monic polynomials can
be defined recursively as $\Phi_1(x) = x-1$ and 
$\D \prod_{i|n} \Phi_i(x) = x^n - 1$
for $n > 1$. The first few  are
easily calculated to be 
$x-1,\  x+1, \ x^2+x+1, \ x^2+1, \dots$\,.
For these and other basic facts, see an algebra text such as \cite{Hung}.

 While it might 
appear that the coefficients of the
cyclotomic polynomials are always $\pm 1$, the presence of 
$2x^7$ in $\Phi_{105}(x)$ shows that this is not invariably the case 
(and indeed is a good counterexample for those students who insist
that the ``law of small numbers\rq\rq\  is universally valid;
%Guy\rq s delightful article 
see \cite{Guy} for further discussion).
% of this and other 
%counter-examples).
Naturally, much work has been done on the values of the coefficients
of $\Phi_n(x)$. 
One amazing fact worthy of 
mention is that
every integer  appears as a coefficient in
some cyclotomic polynomial (see \cite{C}, \cite{Suzuki}).

In this article, we  provide a short and elementary proof 
of the following result:

\

%%%%%%%%%%%%%%%%%%%%%%%%%%%%%%%%%%%%%%%%%%%%%%%%%%%%%%%%%%%%%%%%%%%%%%%%%%%%
\noindent{\bf Theorem 1.}
{\it For $n \geq 3$ the 
middle coefficient of $\Phi_n(x)$
is either zero (when $n$ is a power of $2$) or
an odd integer.}
%%%%%%%%%%%%%%%%%%%%%%%%%%%%%%%%%%%%%%%%%%%%%%%%%%%%%%%%%%%%%%%%%%%%%%%%%%%%

\ 

A similar result 
can be found
 in \cite{LL}, where 
Lam and Leung directly calculate 
the middle coefficient of $\Phi_{pq}(x)$ for distinct 
primes $p$ and $q$ and show it to be
$\pm 1$. This had been done earlier by Beiter \cite{Beiter} for
the case of distinct odd primes. Both papers rely on the partition
of $\phi(pq)/2$ into $rp + sq$. In contrast,
our proof uses 
only some very basic facts about
minimal polynomials.
We also point out that for
$n \not= pq$ the polynomial 
$\Phi_n(x)$ could 
indeed have a middle coefficient different 
from $1$ or $-1$. The first such occurence is at
$n=385$ (giving a middle coefficient of $-3$), 
after which one  sees $5$ at $n=4785$, followed by 
$-7$ at $n=7735$, and 
$19$ at $n=11305$. All these values of $n$ are square-free products of
small odd primes, which is alluded to in \cite{Suzuki}.

Before proceeding with the proof of Theorem 1,
 we do a bit of preliminary work.
The first lemma establishes a useful fact about $\Phi_n(x)$.

\ 

%%%%%%%%%%%%%%%%%%%%%%%%%%%%%%%%
\noindent{\bf Lemma 1.}
{\it If $n \geq 3$ and is odd, then  
$\Phi_n(-1) = 1$.}
%%%%%%%%%%%%%%%%%%%%%%%%%%%%%%%%

\ 

\noindent{\it Proof.}
For $n\geq 3$, 
\[
\prod_{i|n, i>1} \Phi_i(x) =  \frac{x^n-1 }{x-1},
\]
so
(since $n$ odd) 
\[
\prod_{i|n, i>1} \Phi_i(-1) =  \frac{(-1)^n-1 }{(-1)-1} = 1.
\]
Also, $\Phi_3(-1) = 1$. By a simple 
induction argument we  conclude that 
$\Phi_n(-1) = 1$ whenever  $n \geq 3$ and is odd.
\hfill\rule{2mm}{2mm}

\vskip0.2in

Next we review some basic information.
We use
$\zeta_n$
to signify a primitive $n${\rm{th}} root of unity (that is,
$\zeta_n = e^{2\pi ik/n}$ for some $k$ relatively prime to $n$), and 
$f_n(x)$ to denote the minimal polynomial of
$\zeta_n + \zeta_n^{-1}$ (recall that the {\it minimal polynomial\/} of 
an algebraic complex
number $\alpha$ is the
monic polynomial $p(x)$ in ${\Q}[x]$ of smallest degree such that $p(\alpha)=0$). 
It is not hard to show using
elementary methods (see \cite{Lehmer}) that 
$f_n$ has integer coefficients and that when
$n \geq 3$ the degree of $f_n$ 
is
half that of $\Phi_n(x)$. In fact,
\begin{equation}\label{e1}
\qquad \qquad \Phi_n(x) = f_n (x + x^{-1})\cdot x^{\phi(n)/2} \qquad \mbox{($n\geq 3$),} 
\end{equation}
because (after simplifying the right-hand side) the polynomials on both sides of (\ref{e1}) are 
monic, are of degree $\phi(n)$, and have $\zeta_n$ as a root.
The first few such polynomials $f_n$ (for $n \geq 3$) are
easy to derive from (\ref{e1}) and read as follows:

\parbox{3.1cm}{\begin{eqnarray*}
&&f_3(x) = x+1, \\
&&f_4(x) = x, 
\end{eqnarray*}
}
\parbox{4.1cm}{\begin{eqnarray*}
&&f_5(x) = x^2+x-1, \\
&&f_6(x) = x-1,
\end{eqnarray*}
}
\parbox{4cm}{\begin{eqnarray*}
&&f_7(x) = x^3 + x^2 - 2x -1,\\
&&f_8(x) = x^2-2.
\end{eqnarray*}
}

From this, we see that the
constant term in $f_n$ is not always $\pm 1$
(equivalently,  $\zeta_n + \zeta_n^{-1}$ is not
necessarily an {\it algebraic unit}, meaning a unit in the ring of
algebraic integers). However, by
doing a 
careful comparison of the $f_n$ with the Chebyshev polynomials,
Carlitz and Thomas \cite{CT} showed that when $n \geq 3$
and $n$ is not divisible by $4$,
the  constant term in $f_n(x)$ is either $1$ or $-1$.
% (see \cite{CT}), and so in this case
% $\zeta_n + \zeta_n^{-1}$
% is indeed an
% algebraic unit.
For the sake of completeness, we provide  a nonelementary, 
but much shorter, proof of this fact.

\ 

%%%%%%%%%%%%%%%%%%%%%%%%%%%%%%%%
\noindent{\bf Lemma 2.}
{\it If $n\geq 3$ and $n \not\equiv 0\  (\mbox{\rm mod } 4)$,
then 
$\zeta_n + \zeta_n^{-1}$
is an
algebraic unit.}
%%%%%%%%%%%%%%%%%%%%%%%%%%%%%%%%

\ 

\noindent{\it Proof.}
Let $m = n$ for $n$ odd and $m=n/2$ for $n$ even. Note that $m$ 
is itself
odd and $m \geq 3$.
Note as well that 
${\zeta_n}^2$
is a primitive
$m${\rm{th}} root of unity (and thus a root of 
$\Phi_m(x)$). Then
${\zeta_n}^2 + 1$ is a root of 
$\Phi_m(x-1)$, 
which is a monic polynomial  with constant term 
$\Phi_m(-1) = 1$ (by Lemma 1). It follows that
${\zeta_n}^2 + 1$ is 
an algebraic unit, as is $\zeta_n$. Thus,
$\zeta_n + \zeta_n^{-1} = ({\zeta_n}^2 + 1)/\zeta_n$
is likewise
an algebraic unit. 
\hfill\rule{2mm}{2mm}

\vskip0.2in

%%%%%%%%%%%%%%%%%%%%%%%%%%%%%%%%%%%%%%%%%%%%%%%%%%%%%%%%%%%%%%%%%%%%%%%

We are now ready to bring everything together.

\noindent{\it Proof of Theorem 1}.
If
$n = 2^k$, then 
$\Phi_n(x) = x^{2^{k-1}} + 1$,
a polynomial with zero as its
middle coefficient. We proceed assuming
that
 $n$ is not a power of $2$.
 
Note that if $\zeta$ is a primitive $4k${\rm{th}} root of unity, then
$\zeta^2$ is a primitive $2k${\rm{th}} root of unity. Since $\phi(4k) = 2\phi(2k)$, 
we know that 
$\Phi_{4k}(x) = \Phi_{2k}(x^2)$. Since the middle coefficient of
$\Phi_{2k}(x^2)$ is the same as that of
$\Phi_{2k}(x)$, we can further assume without loss of generality 
that $4$ does not divide $n$.

Now letting $f_n(x)$  be the 
 minimal polynomial of
$\zeta_n + \zeta_n^{-1}$, we know 
from Lemma 2 that
$f_n$ has constant coefficient $\pm 1$. Thus, 
we can write 
$f_n(x) =  x^k + a_{k-1}x^{k-1} + \cdots + a_1x \pm 1$
 (for $k=\phi(n)/2$),
and so from 
equation (\ref{e1})
we obtain
\begin{equation}\label{e2}
\Phi_n(x) = 
\left[(x+x^{-1})^k + a_{k-1}(x+x^{-1})^{k-1} + \cdots 
 \pm 1\right]\cdot x^{k}.
\end{equation}
The middle coefficient of $\Phi_n(x)$ is the
coefficient of the $x^k$ term in (\ref{e2}) (recall,
 $k=\phi(n)/2$).
This  number is simply the
 sum of the constant terms appearing in each 
expression
 $a_i(x+x^{-1})^i$ in (\ref{e2}), plus the final $\pm 1$.
  The  constant term in 
 $a_i(x+x^{-1})^i$ is either zero (for $i$ odd) or
$a_i {i \choose {i/2}}$ (for $i$ even). As a result, the middle coefficient 
of $\Phi_n(x)$ is 
\begin{equation}\label{e3}
 \sum_{i = 2j}a_i  {i \choose {i/2}} \pm 1
= 
\sum_{j}a_{2j} { 2j \choose {j}} \pm 1.
\end{equation}
By a familiar identity, 
\[
{ 2j \choose {j}} = 
{ 2j-1 \choose {j-1}} + { 2j-1 \choose {j}} 
= 2 { 2j-1 \choose {j}}.
\] 
Thus
the middle coefficient of
$\Phi_n(x)$ is odd when $n$ is not a power of $2$.
\hfill\rule{2mm}{2mm}

%\vskip0.3in 

\noindent\hrulefill

{
\renewcommand{\baselinestretch}{1.0}
\normalsize

\noindent Dr. Gregory Dresden  \hfill (540) 458-8806

\noindent Department of Mathematics, Robinson Hall \hfill {\tt dresdeng@wlu.edu}

\noindent Washington and Lee University

\noindent Lexington, VA 24450

}

\end{document}